\documentclass[letterpaper, 10 pt, conference]{ieeeconf}  
\usepackage{graphicx} 
\usepackage{fancyhdr}

\IEEEoverridecommandlockouts                              
\title{\LARGE \bf
An enhanced and more realistic tank environment setup for the development of new methods for fish behavioral analysis in aquaculture
}

  \author{Dimitris Voskakis$^{1,\dagger}$, Martin Føre$^{1}$, Eirik Svendsen$^{2}$, Aleksander Perlic Liland$^{2}$, Sonia Rey Planellas $^{3}$,\\ Harkaitz Eguiraun$^{{4},{5}}$, Pascal Klebert $^{2}$ %
\thanks{*This work was supported by the Research Council of Norway through the projects BioWaSys (343277).}
\thanks{$^{1}$Department of Engineering Cybernetics, Faculty of Information Technology and Electrical Engineering, Norwegian University of Science and Technology, Trondheim, Norway}%
\thanks{$^{2}$SINTEF Ocean AS, Trondheim, Norway}%
\thanks{$^{3}$Institute of Aquaculture, School of Natural Sciences, University of Stirling, Stirling FK9 4LA, Scotland, UK}
\thanks{$^{4}$Department of Graphic Design \& Engineering Projects, Faculty of Engineering of Bilbao, University of the Basque Country UPV/EHU}
\thanks{$^{5}$Research Centre for Experimental Marine Biology \& Biotechnology, University of the Basque Country PiE‐UPV/EHU, Plentzia, Spain}
\thanks{$\dagger$Corresponding author: {\tt\small dimitris.voskakis@ntnu.no}}}  

\fancypagestyle{withfooter}{
  
  \fancyfoot[C]{\footnotesize Accepted to the IEEE IROS workshop on Autonomous Robotic Systems in Aquaculture: Research Challenges and Industry Needs}
}

\begin{document}
\maketitle
\thispagestyle{withfooter}
\pagestyle{withfooter}

\begin{abstract}
The aquaculture industry is constantly making efforts to improve fish welfare while maintaining the ethically sustainable farming practises. 
This work presents an enhanced tank environment designed for testing and developing novel combinations of technologies for analyzing and detecting behavioral responses in fish shoals/groups. 
Regular cameras are combined with event cameras and a scanning sonar to comprise a sensor suite that offers a more detailed and complex way of fish observation. 
The modified tank environment is designed to simulate the prevailing conditions on-site at cage based farms, particularly in terms of lighting conditions, while all tank systems and sensors are hidden behind specially designed enclosures, providing a "clean" environment (open arena) less likely to impact the fish behavior. 
The proposed sensor suite will be tested and demonstrated in the modified tank environment to benchmark its ability in monitoring fish, after which it will be adapted for use in a more industrially relevant situation with open cages.

\noindent\textbf{\emph{Keywords:}} \textnormal{\emph{Aquaculture}, \emph{Fish behavior}, \emph{Light Conditions}, \emph{Precision Fish Farming (PFF)}, \emph{Computer vision, Fish Welfare}}

\end{abstract}

\section{Introduction}
\subsection{Motivation}
Over the last decades, the Norwegian aquaculture industry has experienced significant technological developments and breakthroughs, cementing its leading global role within intensive fish production as a sustainable food source for human consumption \cite{paisley_overview_2010, afewerki_innovation_2023}.
Atlantic salmon is the most commercially valuable species farmed in Norway, and the industry has made significant efforts to enhance fish welfare, with the dual purpose of maintaining and improving the quality of the final product and ensuring ethical production practices. 
Recent trends have also sought to improve production by replacing labour-intensive procedures and practices with new technologies, autonomy and more objective methods, which is in tune with precision farming methods \cite{fore_precision_2018}.
Systems utilizing hydroacoustics, cameras and biosensors have been key elements in this transition, offering valuable insights into fish biology (such as behavior, morphometric characteristics) and the dynamics of the production environment \cite{voskakis_deep_2021, sadoul_new_2014, brijs_bio-sensing_2021, kristmundsson_fish_2023, saberioon_application_2017}.
    
Several studies have aspired to estimate the behavioral responses \cite{georgopoulou_tracking_2021, eguiraun_entropy_2023} and morphometric characteristics \cite{voskakis_deep_2021, monkman_accurate_2020, tseng_automatic_2020} of fish, many of which are based on artificial vision.  
This approach is not limited to standalone procedures, and optical technologies have also been mounted on robotic systems and used for underwater operations in aquaculture cages, thereby enabling their use as perception tools for navigation and obstacle avoidance \cite{kelasidi_robotics_2022,skaldebo_modeling_2023}. 
The camera technology required to set up machine vision applications are today affordable and robust enough for subsurface applications. 
However, in most cases, the cameras must be calibrated to achieve precise estimations, resulting in a set of optimal parameters to use in further processing \cite{voskakis_design_nodate}, \cite{zhang_flexible_2000} \cite{shortis_calibration_2015}, \cite{lavest_dry_2003}. 
After calibration, the images are typically subjected to one or more of the increasing number of AI-based methods for processing camera streams to yield the desired output.
Using camera based methods often also require  artificial lights, since the object of interest needs to be illuminated. 
Such lights must however be applied with care, as submerged lighting is used as an active measure in the aquaculture industry due to its impacts on fish growth, welfare, smoltification and changes in maturation stages \cite{hansen_effects_2017}, \cite{endal_effects_2000}.

\subsection{Previous Works}
\subsubsection{Conventional cameras}
In \cite{voskakis_deep_2021, shi_automatic_2020, meng-che_chuang_tracking_2015} conventional vision cameras were used to analyze the morphometric parameters of the fish by submerging a waterproof camera in open sea cages. 
To acquire accurate data using such methods, a calibration process is required. 
In \cite{voskakis_design_nodate, yau_underwater_2013, voskakis2024modeling,kwon_applicability_nodate, she_adjustment_2019} various calibration methods highlight the need and process for optimal extraction of camera parameters to obtain accurate estimates.

\subsubsection{Hydroacoustic methods}
Hydroacoustic technologies are considered effective tools for observing fish in their natural habitat or during captivity in the aquaculture \cite{fore_precision_2018}. 
Research studies have explored the potential of such systems, highlighting several useful results \cite{kristmundsson_fish_2023}, \cite{tao_hydroacoustic_2010}, \cite{knudsen_hydroacoustic_2004}. 
Most of these previous efforts using hydroacoustic methods have been done in open sea-cages, and such tools may often be unsuitable for tank deployment due to multipathing issues. 
However, some light-weight scanning devices have recently shown promise in data collection at short ranges \cite{zhang2024farmed}, implying their potential also in more confined spaces.


\subsubsection{Artificial lights}
Since optical cameras and hydroacoustic devices all have their advantages and drawbacks, a combination of these technologies would result in a robust and multi-modal system for monitoring fish.
This approach could leverage the strengths of the aforementioned technologies , thereby achieving both higher accuracy in predictions, and improved real-time performance.
Since the system will feature optical methods, appropriate ambient lighting must be maintained during data collection. 
However, illumination must be set up properly to avoid compromising production since lighting may affect maturity, development and growth  \cite{endal_effects_2000}, \cite{handeland_low_2013}, \cite{imsland_effect_2017}. 

\subsubsection{Project aims}
Our main aim is to establish non-invasive sensing methods for assessing stress responses on a shoal/group level for fish monitoring and behavioral responses in aquaculture settings.
Using the fish as "sensors", we will thus provide a foundation for enabling warning systems that inform about changes in the fish state as early as possible.  
This will be done through a series of experiments where the fish are exposed to different external stressors (e.g., environmental changes, flashing light).
During the exposure, fish behaviour will be monitored using different sensor types to capture several different modes of the response patterns. 
The first trials are aimed at identifying the sensors best able to capture elicited stress responses. 
These trials will be conducted in an enhanced tank environment that is designed to simulate a sea-cage replica by mimicking the prevailing light, where salmon will encounter in commercial aquaculture cages. 
When a sensor setup suitable for capturing group level stress responses in salmon is identified, followup experiments will seek to apply the same sensor suite in commercial cages, thereby validating the setup for industrial use.

\section{Materials and Methods}
\subsection{Tank setup and design}
Several modifications will be made to create a tank environment that simulates the prevailing conditions in commercial cages, e.g., by mimicking the subsurface light levels. 
An octagon shape tank with internal dimensions of 2m x 2m x1m (apx 4m$^3$) will be used as a basis for  building the setup. To modify light conditions, a flexible submersible LED line (Table \ref{tab:light_specs}) will be placed close to the surface facing downwards, illuminating the tank with artificial daylight (5600K to 6500K Color Temperature). 
Homogeneous illumination inside the tank was achieved by bending the LED line horizontally, and attaching it to a support frame mounted along the inner perimeter of the tank, and using a relatively wide opening angle (112.3$^\circ$). 
The support frame is used to stabilise the vertical position of the LED line, and is equipped with vertical screws that enable depth adjustments according to the desired experimental setup.
Furthermore, the LEDs were set up with a color rendering index (CRI) of 80 as this would result in underwater visibility conditions resembling those with natural light in sea-cages.

\begin{table}[htbp]
    \centering
    \caption{Technical specifications of the submersible LED lights.}
    \scalebox{0.9}{
    \begin{tabular}{|c|c|}
        \hline
        \textbf{Parameters} & \textbf{Values}\\
        \hline
        Color temperature & 5700 K\\
        Power supply & 24V DC (15 W)\\
        CRI & 80\\
        Average beam angle & 112.3$^\circ$\\
        Lumen/m & 700\\
        \hline
    \end{tabular}
    } \vspace{-5pt}
    \label{tab:light_specs}
\end{table}

\begin{figure}[!h]
    \centering
    \includegraphics[width=0.8\linewidth]{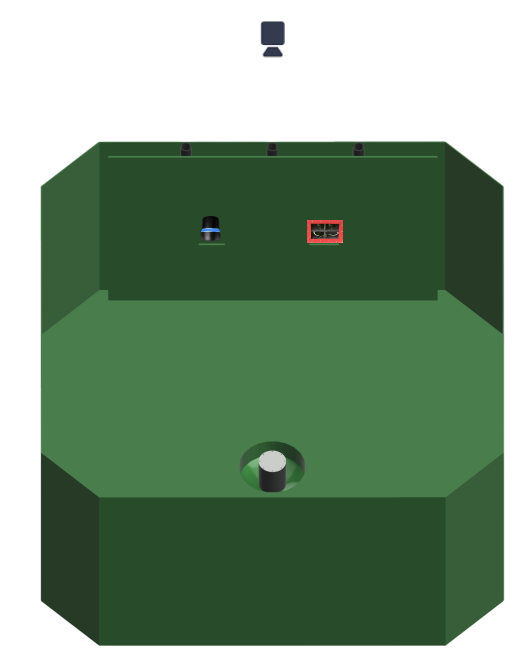}
    \caption{3D demonstration of the tank suite.}
    \label{fig:3d_tank}
\end{figure}

In addition to simulating on-site light conditions, the tank was designed to present the fish within an open arena  where the various instruments/sensors, water inlets and other objects present in the tank, were not visible for the fish. 
The idea behind this was to both reduce the risk of the fish colliding with the components inside the tank, and to inhibit the potential impact of the presence of these on fish responses. 
This was realised by installing a camouflage frame that hides the inlets/outlets of the tank (Fig.\ref{fig:3d_tank}).
The camouflage frame was used to conceal two of the instruments (Stereo video - Right and Sonar - Left) used to monitor the fish, ensuring high-quality data acquisition without disturbing the fish during the experiments. 
A hole at the bottom of the tank will be used to deploy the last instrument (event camera), thereby hiding also this from the fish. Subsequently, a monocular camera wiil be placed on the ceiling, monitoring the fish behavior from a different point of view.
In sum, this resulted in the setup better emulating the environment in a commercial cage, since farmers try to prevent the presence of external elements inside the production volume as well as possible.  

\subsection{Sensing technologies}
Four instruments will be used to capture and analyze the behavior of the fish, featuring three types of cameras and one sonar. 
A stereo and a monocular camera will be used to acquire high resolution picture frames and detect the potential biological responses and behavior of the fish, using conventional AI-based processing methods. 
Since shoal level responses are the main focus of this work, the outputs from the stereo video system will mainly be aggregated over time to achieve values representing groups rather than individuals. 
Properties of interest here would be external morphometrics, swimming speeds and accelerations, and positioning/distribution in the tank volume. 

The other optical method used in the setup is a newly emerged type of camera called an event camera or Dynamic Vision Sensor (DVS). It is a new technology originally designed for the automotive industry, which we believe it may have many applications within the aquaculture. 
Event cameras are mostly used to detect the immediate environment around aerial vehicles (such as drones) and cars \cite{gallego_event-based_2022}, \cite{gallego2017accurate}, \cite{rebecq2019high}, and thus improve their perception and hence reaction speed. 
This camera type differs from conventional cameras (such as those used in the stereo video system here) in that it can output pixel-level brightness changes instead of capturing regular RGB frames.
The advantage of this is that the camera will only output data when there exist changes in the captured image, but only for the pixels that are changed. 
This enables the camera to report data at a far higher speed than possible through post-processed RGB frames. 
In addition to improving the response time, this also contributes to effectually inhibit features that may compromise outputs from conventional processing such as motion blur.   
However, to our knowledge, event cameras have never before been applied in aquaculture settings, making this application a world first. 

\begin{figure}[h!]
    \centering
    \includegraphics[width=0.7\linewidth]{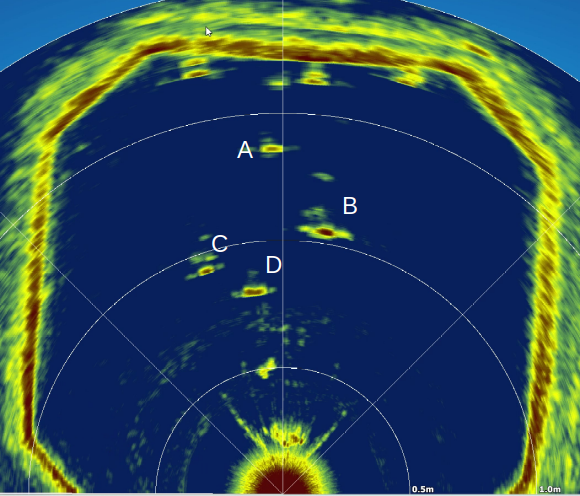}
    \caption{Deployment of a Ping360 - Dummy objects (A, B,C, D).}
    \label{fig:Ping360}
\end{figure}

To also explore the potential of acoustic methods, the last monitoring tool used in the sensor suite is a Ping360 sonar  (BlueRobotics Inc.).
Previous studies have shown this 360 $^\circ$ scanning sonar to be an efficient sensor for detecting fish in aquaculture settings \cite{zhang2024farmed}.
However, multipathing has previously been found to be a major issue when using acoustic devices in tanks since the walls tend to have high target strength. 
A preliminary experimental work using dummy objects in the tank (i.e., inflated balloons moored to metal objects at the bottom) indicated the potential use of this system. 
These tests proved the potential to capture both the circumference of the tank and the objects placed in the tank volume, implying it should work as intended in the experiment (Fig.\ref{fig:Ping360}). Underwater acoustic imaging could be a new way to visualise and monitor fish behaviour.



\section{Conclusions}
This work explores the development of an enhanced tank environment that mimics the environmental conditions, particularly light, prevailing in commercial sea cages to enable tank experiments in more realistic conditions. 
Furthermore, a combination of systems will be tested in their ability to analyze the behavioral patterns of salmon exposed to various stress factors. 
The sensor suite used for this features stereo cameras, an event camera and a sweeping 360$^\circ$ sonar.
The ultimate goal of the work is to identify the sensor combinations providing best insight into behavioral responses in fish, and then use these to test and quantify responses in salmon exposed to various stressors. 
A follow up experiment will build on the systems, observations and knowledge obtained in the tank experiments to design a similar setup for use in commercial cage environments. 
Joint analyzes will then be used to assess the estimated parameters from both environments and find similarities that could pave the way for new technological solutions for shoal level monitoring in aquaculture.

\section*{Acknowledgement}
This work was supported by the Research Council of Norway through the projects BioWaSys (343277).

\bibliographystyle{ieeetr}
\bibliography{Library.bib}

\end{document}